\documentclass[12pt,a4paper]{article}
\usepackage[latin2]{inputenc}
\usepackage{amsmath}
\usepackage{amsfonts}
\usepackage{amssymb}
\usepackage{graphicx}
\usepackage[left=2.00cm, right=2.00cm, top=2.50cm, bottom=2.50cm]{geometry}
\usepackage[T1]{fontenc}

\def\Q{{\mathbb Q}}
\def\Z{{\mathbb Z}}

\def\O{{\cal O}}

\newtheorem{lemma}{Lemma}

\newtheorem{corollary}[lemma]{Corollary}

\title{
Power integral bases\\in cubic and quartic extensions of\\real quadratic fields
}
\author{
Istv\'{a}n Ga\'{a}l\thanks{
        Research supported in part by K115479 from the
        Hungarian National Foundation for Scientific Research				
				and by the EFOP-3.6.1-16-2016-00022 project. 
				The project is co-financed by the European Union and the European Social Fund.	
                         },\; \\
{\small University of Debrecen, Mathematical Institute} \\
{\small H--4002 Debrecen Pf.400., Hungary,} 
{\small e--mail: gaal.istvan@unideb.hu},
 \\ \\
L\'aszl\'o Remete\thanks{
        Research supported by the \'UNKP-18-3 new national excellence program of the ministry of human capacities.}\; \\
{\small University of Debrecen, Mathematical Institute} \\
{\small H--4002 Debrecen Pf.400., Hungary,} 
{\small e--mail: remete.laszlo@science.unideb.hu}
}

\begin{document}
\baselineskip=17pt

\maketitle
\thispagestyle{empty}

\renewcommand{\thefootnote}{\arabic{footnote}}
\setcounter{footnote}{0}


\begin{abstract}
Investigations of monogenity and power integral bases 
were recently extended from the absolute case (over $\Q$) to the relative case
(over algebraic number fields). 
Formerly, in the relative case we only succeeded to calculate generators
of power integral bases when the ground field is an
imaginary quadratic field. 
This is the first case when we consider monogenity in the more difficult case,
in extensions of real quadratic fields.
We give efficient algorithms for calculating generators of power integral 
bases in cubic and quartic extensions of real quadratic fields, more exactly in
composites of cubic and quartic fields with real quadratic fields. In case
the quartic field is totally complex, we present an especially simple
algorithm.

We illustrate our method with two detailed examples.
\end{abstract}

\newpage

\section{Introduction}
Let $K$ be an algebraic number field with ring of integers $\Z_K$.
This field is {\it monogenic} 
if $\Z_K$ is a simple ring extension of $\Z$, that is
$\Z_K=\Z[\vartheta]$ with some $\vartheta\in\Z_K$. In this case
$\{1,\vartheta,\ldots,\vartheta^{n-1}\}$ is an integral basis of $K$, called
{\it power integral basis}.
Monogenity of number fields and the calculation of 
generators of power integral bases
is a classical topic of algebraic number theory.

In several cases, even in higher degree fields there are methods to
exclude monogenity of the field, cf.
\cite{gop}, \cite{olaj}, \cite{composite}, \cite{simlestsextic}.
If these methods are applicable, then they are very easy to use.

It is much more complicated to explicitly calculate the 
generators of power integral bases c.f. \cite{book}. 
We have efficient algorithms for cubic and quartic number fields, \cite{gsch},
\cite{gppsim}, general but less efficient methods 
for quintic and sextic fields \cite{s5}, \cite{s6}. 
We only have partial results for higher degree fields \cite{compos}, 
\cite{relquartic}, \cite{degnine}.

We also considered monogenity and power integral bases in the 
{\it relative case} 
\cite{relcubic}, \cite{relquartic}, \cite{gsz2}, \cite{jalfa}, \cite{fj18}.
The element $\vartheta$ generates a {\it relative power integral basis} of $K$ over the subfield $M$ if
$\Z_K=\Z_M[\vartheta]$ with some $\vartheta\in\Z_K$.
Until now we only succeeded to consider number fields and 
parametric families of number fields that are relative extensions of
{\it imaginary} quadratic fields. In this paper we extend these results to
the case when the ground field is a {\it real} quadratic number field.

We present efficient algorithms for determining generators of power
integral bases in sextic and octic fields that are composites
of real quadratic fields with cubic and quartic fields, respectively.
The algorithm reduces the problem to solving relative Thue equations
and finding roots of polynomials. In the octic case, if the quartic
field involved is totally complex, then we considerably simplify 
the calculation.

We illustrate our method with two detailed examples.

\section{Preliminaries}
\label{prelim}

Herewith we recall some results of \cite{jalfa} that we shall use in the
following.

Let $M$ be an algebraic number field of degree $m$ and $K$ an extension of
$M$ with $[K:M]=k$, with rings of integers $\Z_M$ and $\Z_K$, respectively.
We have $[K:\Q]=k\cdot m$.
We assume that there exist a relative integral basis of $\Z_K$ over $\Z_M$.
(As we shall see the existence of a power integral 
basis of $\Z_K$
implies the existence of a relative power integral basis.)

Denote by $D_K$ and $D_M$ the discriminants of $K$ and $M$, respectively.
The {\it index} of a primitive element $\alpha$ of $\Z_K$ is
\begin{equation}
I(\alpha)=(\Z_K^+ : \Z[\alpha]^+)=\sqrt{\left|\frac{D(\alpha)}{D_K}\right|},
\label{absind}
\end{equation}
and we also have
\begin{equation}
I(\alpha)=(\Z_K^+ : \Z_M[\alpha]^+)\cdot 
(\Z_M[\alpha]^+: \Z[\alpha]^+),
\label{ind}
\end{equation}
where the indices of the additive groups of the corresponding rings are calculated.
The first factor is the {\it relative index} of $\alpha$: 
\[
I_{K/M}(\alpha)=(\Z_K^+ : \Z_M[\alpha]^+).
\]
Denote by $D_{K/M}$ the relative discriminant of $K$ over $M$.
As it is well known
\begin{equation}
D_K=N_{M/\Q}(D_{K/M})\cdot D_M^{[K:M]}.
\label{discr}
\end{equation}
Denote by $\gamma^{(i)}$ the conjugates of any $\gamma\in M$
($i=1,\ldots,m$). Let $\delta^{(i,j)}$ 
be the images of any $\delta\in K$ under the automorphisms of $K$ leaving the
conjugate field $M^{(i)}$ element-wise fixed ($j=1,\ldots,k$).
Then for any primitive element $\alpha\in \Z_K$ we have
\begin{equation}
I_{K/M}(\alpha)=
\frac{\sqrt{|N_{M/\Q}(D_{K/M}(\alpha))|}}{\sqrt{|N_{M/\Q}(D_{K/M})|}}
=\frac{1}{\sqrt{|N_{M/\Q}(D_{K/M})|}}
\cdot
\prod_{i=1}^m\;\;\prod_{1\leq j_1<j_2\leq k}\left|\alpha^{(i,j_1)}-\alpha^{(i,j_2)}\right|.
\label{relind}
\end{equation}
Further, by (\ref{absind}), (\ref{ind}), (\ref{discr})  and (\ref{relind})
we have
\begin{equation}
J(\alpha)=(\Z_M[\alpha]^+: \Z[\alpha]^+)=
\frac{1}{\sqrt{|D_M|}^{\;[K:M]}}
\cdot
\prod_{1\leq i_1<i_2\leq m}\;\;\prod_{j_1=1}^k \prod_{j_2=1}^k
\left|\alpha^{(i_1,j_1)}-\alpha^{(i_2,j_2)}\right|.
\label{ind2}
\end{equation}
Equation (\ref{ind}) implies that 
a primitive element $\alpha\in\Z_K$ generates a power integral basis of $\Z_K$, 
if and only if 
\begin{equation}
I_{K/M}(\alpha)=1, \;\; {\rm and}\;\;
J(\alpha)=(\Z_M[\alpha]^+: \Z[\alpha]^+)=1.
\label{pp}
\end{equation}
Therefore
if $\alpha$ generates a power integral basis of $\Z_K$, then it
generates a relative power integral basis of $\Z_K$ over $M$.

It is well known that generators or relative power
integral bases are determined up equivalence, that is up 
to multiplication by a unit in $M$
and up to translation by element of $\Z_M$. More exactly,
if $\alpha$ generates a power integral basis of $\Z_K$, then
\[
\alpha=A+\varepsilon \cdot \alpha_0,
\]
where $\alpha_0$ is a generator of a relative power integral basis of $\Z_K$ over $M$,
$\varepsilon$ is a unit in $M$ and $A\in\Z_M$.

Summarizing, in order to determine all generators of power integral bases of
$\Z_K$ we have to determine up to equivalence all generators $\alpha_0\in\Z_K$ of relative power integral bases of $\Z_K$ over $M$. Then for each $\alpha_0$ 
we can use $J(\alpha)=1$ to determine
the unit $\varepsilon\in M$ and $A\in\Z_M$ so that $I(\alpha)=1$ be satisfied.
\\

Let $\{b_1=1,b_2,\ldots,b_m\}$  be an integral basis of $M$. Then the above $A$
can be represented in the form 
$
A=a_1+a_2b_2+\ldots +a_m b_m.
$
The index of $\alpha$ is independent from $a_1$,
hence, using $J(\alpha)=1$, an equation
of degree $k^2m(m-1)/2$, we have to determine $\varepsilon$ and $a_2,\ldots,a_k$.
This task can became very complicated. However if $M$ is an imaginary
quadratic field, then there are only finitely many units $\varepsilon$ in $M$ and
we get a polynomial equation in one variable $a_2$. 
This is what we have used in some of our recent results e.g. \cite{s62}.
In the present paper we consider the more difficult case when
$M$ is a real quadratic field in which case we have two unknowns,
$\varepsilon$ and $a_2$. We can overcome this difficulty, if $J(\alpha)=1$
factorizes into two factors, to obtain two equations instead of
just one.

\section{Factors of index forms of composite fields}

Let $L$ and $M$ be number fields of degree $\ell=[L:\Q]$ and $m=[M:\Q]$,
respectively. Let $K$ be the composite of $L$ and $M$,
that is $K=LM$. 
Denote by $\Z_K,\Z_L,\Z_M$ the rings of integers of $K,L,M$, respectively.
Assume $[K:\Q]=\ell\cdot m$.
Denote by $\vartheta_L,\vartheta_M$ generator elements of 
$L,M$, respectively. For any $\vartheta\in K$ denote by $\vartheta^{(i,j)}$
the conjugate of $\vartheta$ corresponding to the conjugates $\vartheta_M^{(i)}$
and $\vartheta_L^{(j)}$, for $1\leq i\leq m$ and $1\leq j\leq \ell$.

For any primitive element $\alpha\in\Z_K$ we have
\begin{equation}
I(\alpha)=I_{K/M}(\alpha) J_M(\alpha)
\label{iM}
\end{equation}
where
\begin{equation}
I_{K/M}(\alpha)=\frac{1}{\sqrt{|N_{M/\Q}(D_{K/M})|}}
\cdot
\prod_{i=1}^m\;\;\prod_{1\leq j_1<j_2\leq \ell}
\left|\alpha^{(i,j_1)}-\alpha^{(i,j_2)}\right|.
\label{relM}
\end{equation}
and
\begin{equation}
J_M(\alpha)=\frac{1}{\sqrt{|D_M|}^{[K:M]}}
\cdot
\prod_{1\leq i_1<i_2\leq m}\;\;\prod_{j_1=1}^{\ell} \prod_{j_2=1}^{\ell}
\left|\alpha^{(i_1,j_1)}-\alpha^{(i_2,j_2)}\right|.
\label{JM}
\end{equation}

On the other hand 
\begin{equation}
I(\alpha)=I_{K/L}(\alpha) J_L(\alpha)
\label{iL}
\end{equation}
where
\begin{equation}
I_{K/L}(\alpha)=\frac{1}{\sqrt{|N_{L/\Q}(D_{K/L})|}}
\cdot
\prod_{j=1}^{\ell}\;\;\prod_{1\leq i_1<i_2\leq m}
\left|\alpha^{(i_1,j)}-\alpha^{(i_2,j)}\right|.
\label{relL}
\end{equation}
and
\begin{equation}
J_L(\alpha)=\frac{1}{\sqrt{|D_L|}^{[K:L]}}
\cdot
\prod_{1\leq j_1<j_2\leq \ell}\;\;\prod_{i_1=1}^{m} \prod_{i_2=1}^{m}
\left|\alpha^{(i_1,j_1)}-\alpha^{(i_2,j_2)}\right|.
\label{JL}
\end{equation}
Further, set 
\[
D_{L,M}=\frac{D_K}{  N_{M/\Q}(D_{K/M})  \cdot N_{L/\Q}(D_{K/L}) }
\]
and
\begin{equation}
J_{L,M}(\alpha)=\frac{1}{\sqrt{|D_{L,M}|}}
\prod_{(i_1,j_1,i_2,j_2)\in H}
\left|\alpha^{(i_1,j_1)}-\alpha^{(i_2,j_2)}\right|
\label{}
\end{equation}
where 
$H=\{(i_1,j_1,i_2,j_2)\; | \;1\leq i_1<i_2\leq m, 1\leq j_1,j_2\leq \ell, j_1\neq j_2\}$.

Obviously we have
\begin{lemma}
\label{ifactor}
\[
I(\alpha)=I_{K/M}(\alpha) \cdot I_{K/L}(\alpha) \cdot J_{L,M}(\alpha).
\]
\end{lemma}

\begin{corollary}
\label{ifactorcorr}
If 
\[
I(\alpha)=1,
\]
then
\[
I_{K/M}(\alpha)=1,\;\; I_{K/L}(\alpha)=1,\;\;  J_{L,M}(\alpha)=1.
\]
\end{corollary}

\noindent
{\bf Proof of Corollary \ref{ifactorcorr}}.
If $I(\alpha)=1$, then by (\ref{pp}) we have 
$I_{K/M}(\alpha)=1$ and $I_{K/L}(\alpha)=1$ which implies
by Lemma \ref{ifactor} that
\[
J_{L,M}(\alpha)=\frac{I(\alpha)}{I_{K/M}(\alpha)\cdot I_{K/L}(\alpha)}=1.
\]
$\Box$.

To the indices and relative indices there correspond index forms and relative index forms. Let 
$\{b_1=1, b_2,\ldots, b_{\ell}\}$ be linearly independent element of $\Z_L$,
let
$\{f_1=1, f_2,\ldots, f_m\}$ be linearly independent element of $\Z_M$.
Then any $\alpha\in\Z_K$ can be written in the form
\begin{equation}
\alpha=\frac{1}{d}\sum_{r=1}^m\sum_{s=1}^{\ell}x_{rs}b_rf_s,
\label{repr}
\end{equation}
where $x_{rs}\in\Z, (1\leq r\leq m,1\leq s\leq \ell)$ and $0\ne d\in\Z$
is a common denominator. Then
\[
\alpha^{(i,j)}=\frac{1}{d}\sum_{r=1}^m\sum_{s=1}^{\ell}x_{rs}b_r^{(i)}f_s^{(j)}
\]
for $1\leq i\leq m,1\leq j\leq \ell$.
Keeping the $x_{rs}$ as variables and calculating the indices and 
relative indices with this representation 
we obtain polynomials in $x_{rs}$. For brevity set
$\underline{x}=(x_{11},\ldots,x_{1\ell},\ldots,x_{m1},\ldots,x_{m\ell})$.
Define 
$I(\underline{x}),
I_{K/M}(\underline{x}),
I_{K/L}(\underline{x})$ and
$J_{L,M}(\underline{x})$ by
\[
I(\alpha)=\frac{1}{d^{m\ell(m\ell-1)/2}}\; \cdot \; I(\underline{x}),
\]
\[
I_{K/M}(\alpha)=\frac{1}{d^{m\ell(\ell-1)/2}}\; \cdot \; I_{K/M}(\underline{x}),
\]
\[
I_{K/L}(\alpha)=\frac{1}{d^{m\ell(m-1)/2}}\; \cdot \; I_{K/L}(\underline{x}),
\]
\[
J_{L,M}(\alpha)=\frac{1}{d^{(m^2\ell^2-m^2\ell-m\ell^2+m\ell)/2}}\; \cdot \; 
J_{L,M}(\underline{x}).
\]
It is known that $I(\underline{x})$,
$I_{K/M}(\underline{x})$ and $I_{K/L}(\underline{x})$
have integer coefficients as polynomials in 
$x_{11},\ldots,x_{1\ell}$, $\ldots$, $x_{m1},\ldots,x_{m\ell}$.

This factorization of the index form is what we shall utilize.
We consider composites of cubic and quartic fields with real quadratic fields.
If we can determine the elements of relative index 1 over $M$,
then we have two more equations to calculate the two unknowns
$\varepsilon$ and $a_2$ (c.f. end of Section \ref{prelim}),
to determine the elements of absolute index 1.

\section{Composites of cubic fields with real quadratic fields}
\label{63}

Let $M$ be a real quadratic number field with integral basis $\{1,\omega\}$
and fundamental unit $\varepsilon>1$.
Let $L=\Q(\xi)$ be a cubic number field, with a generator element $\xi\in\Z_L$.
In this case the representation (\ref{repr}) is written for simplicity in the
form
\begin{equation}
\alpha=\frac{x_0+x_1\xi+x_2\xi^2+y_0\omega+y_1\omega \xi+y_2\omega \xi^2}{d},
\label{repr6}
\end{equation}
where $x_i,y_i\in\Z\; (0\leq i\leq 2)$ and $0<d\in\Z$ is a common denominator.
For $1\leq i\leq 2$, \mbox{$1\leq j_1<j_2\leq3$} we have
\[
\alpha^{(i,j_1)}-\alpha^{(i,j_2)}=
\frac{1}{d}
(\xi^{(j_1)}-\xi^{(j_2)})
\left[(x_1+\omega{^{(i)}}y_1)+(\xi^{(j_1)}+\xi^{(j_2)})(x_2+\omega{^{(i)}}y_2)\right].
\]
Let $a\in\Z$ denote the coefficient of the quadratic term in the
cubic minimal polynomial of $\xi$. Then
\[
I_{K/M}(\alpha)=\frac{1}{d^{6}} \cdot 
\frac{1}{\sqrt{|N_{M/\Q}(D_{K/M})|}}
\cdot
\prod_{i=1}^2\;\;\prod_{1\leq j_1<j_2\leq 3}
\left|\alpha^{(i,j_1)}-\alpha^{(i,j_2)}\right|=
\]
\begin{equation}
\frac{1}{d^{6}} \cdot 
\frac{|D(\xi)|}{   \sqrt{|N_{M/\Q}(D_{K/M})|}}\; \;
|N_{M/\Q}\left(N_{K/M}(X-(a+\xi)Y)\right)|,
\label{relMx}
\end{equation}
with integer variables $X=x_1+\omega y_1,Y=x_2+\omega y_2$ in $M$.
If $\alpha$ generates a power integral basis in $K$, then by 
Corollary \ref{ifactorcorr} we have $I_{K/M}(\alpha)=1$. In view of 
(\ref{relMx}) this gives rise to a {\it relative Thue equation:}
\begin{equation}
N_{M/\Q}\left(N_{K/M}(X-(a+\xi)Y)\right)=\pm m_0
\label{r3}
\end{equation}
with 
\[
m_0=\frac{d^6 \sqrt{|N_{M/\Q}(D_{K/M})|}}{|D(\xi)| }
\]
which is an integer, since the norm of an algebraic integer is on the 
left side.

There is an efficient algorithm to calculate solutions of relative Thue
 equations, see \cite{gp}, \cite{book}. The procedure involves the application
of Baker's method, reduction algorithms and enumeration of small exponents.
Note that there is also a fast algorithm \cite{relthuesmall}
for calculating the "small" solutions (of size, say $<10^{200}$)
of relative Thue equations.

Equation (\ref{r3}) enables us to determine the solutions $X,Y\in\Z_M$
up to a unit factor, that is to determine all possible 
$X_0,Y_0\in\Z_M$
such that 
\begin{equation}
X=x_1+\omega y_1=\pm \varepsilon^h X_0,\;\; 
Y=x_2+\omega y_2=\pm \varepsilon^h Y_0
\label{yyyy}
\end{equation}
with some $h\in\Z$.
Therefore we obtain 
\begin{eqnarray}
x_1&=& \pm 
\frac{\varepsilon^h \overline{\omega}X_0
-\overline{\varepsilon}^h \omega \overline{X_0}}
{\overline{\omega}-\omega},\nonumber\\
y_1&=&\pm
\frac{\varepsilon^h X_0
-\overline{\varepsilon}^h  \overline{X_0}}
{\omega-\overline{\omega}},\label{xy}\\
x_2&=& \pm 
\frac{\varepsilon^h \overline{\omega}Y_0
-\overline{\varepsilon}^h \omega \overline{Y_0}}
{\overline{\omega}-\omega},\nonumber\\
y_2&=&\pm
\frac{\varepsilon^h Y_0
-\overline{\varepsilon}^h  \overline{Y_0}}
{\omega-\overline{\omega}},\nonumber
\end{eqnarray}
where overline denotes the conjugate of an element of $M$.
We set $e=\varepsilon^h$. Then $\overline{\varepsilon}^h=\pm e^{-1}$
the sign depending on the norm of $\varepsilon$.

We have two additional equations:
\[
I_{K/L}(\alpha)=1,\;\; J_{L,M}(\alpha)=1,
\]
whence
\[
I_{K/L}(\underline{x})=I_{K/L}(x_1,x_2,y_0,y_1,y_2)=d^3,\;
J_{L,M}(\underline{x})=J_{L,M}(x_1,x_2,y_0,y_1,y_2)=d^6.
\]
If we substitute the expressions (\ref{xy}) of $x_1,y_1,x_2,y_2$ 
into these
equations, 
and multiply them by $e^3$ and $e^6$, respectively, 
then we obtain polynomial equations in $e,y_0$:
\begin{equation}
F_2(e,y_0)=0,\; F_3(e,y_0)=0.
\label{ee}
\end{equation}
Here $F_2$ is of degree 6 in $e$ and of degree 3 in $y_0$, 
$F_3$ is of degree 12 in $e$ and of degree 6 in $y_0$.
Let $F_4(y_0)$ be the resultant of $F_2$ and $F_3$ with respect to $e$.
This is a polynomial equation in $y_0$ with rational coefficients, of degree 72. 
It is easy to find the integer roots $y_0$ of $F_4(y_0)=0$.
For a given $y_0$ we calculate the real roots of
\[
F_2(e,y_0)=0
\]
in $e$ and for all real roots $e=\gamma$ we check if
\[
\frac{\ln(|\gamma|)}{\ln\varepsilon}
\]
is an integer. If so, then this value is a candidate for $h$.
We obtain $x_1,y_1,x_2,y_2$ from (\ref{xy}).
We have to check if $\alpha$ of (\ref{repr6}) is an algebraic
integer with some $x_0\in\Z$
and if it has $I(\alpha)=1$.

\section{Example 1}

Let $M=\Q(\sqrt{19})$ with fundamental unit $\varepsilon=170+39\sqrt{19}$.
We set $\omega=\sqrt{19}$.
Let $L=\Q(\xi)$ where $\xi$ has minimal polynomial $x^3-3x+17$.
The composite field $K=LM=\Q(\xi,\sqrt{19})$ is generated
by $\vartheta=\xi\sqrt{19}$ 
with minimal polynomial $x^6-114x^4+3249x^2-1982251$.
Considering the integral basis of $K$ we may represent all $\alpha\in\Z_K$
in the from
\begin{equation}
\alpha=\frac{x_0+x_1\xi+x_2\xi^2+y_0\omega+y_1\omega \xi+y_2\omega \xi^2}{19}
\label{ggg3}
\end{equation}
with $x_i,y_i\in\Z\ (0\leq i\leq 2)$.
Our purpose is to determine
all $\alpha\in\Z_K$ having index 1, that is, generating a power integral basis.

We have 
\[
D_M=4\cdot 19,\; D_L=-3^4\cdot 5\cdot 19,\; D_K=2^6\cdot 3^8\cdot 5^2\cdot 19^3.
\]
By 
\[
D_K=N_{M/\Q}(D_{K/M})\cdot D_M^3,\;\; D_K=N_{L/\Q}(D_{K/L})\cdot D_L^2
\]
we obtain
\[
N_{M/\Q}(D_{K/M})=3^8\cdot 5^2, \;\; N_{L/\Q}(D_{K/L})=2^6\cdot 19.
\]
The relative Thue equation (\ref{r3}) has the form
\begin{equation}
N_{M/\Q}(N_{K/M}(X-\xi Y))=\pm 19^5.
\label{rt}
\end{equation}
Up to associates there are three elements of norm $\pm 19^5$ in $\Z_K$:
\[
\gamma_1=\omega+\omega\xi,\; 
\gamma_2=2\omega-\omega\xi,\; 
\gamma_3=\frac{-77\omega+28\omega\xi-9\omega\xi^2}{19}.
\]
The fundamental units of $K$ are
\[
\eta_1=\varepsilon,\; \eta_2=-6+3\xi-\xi^2
\]
and
\[
\eta_3=138357-88194\xi-45680\xi^2+\frac{603932}{19}\omega-
\frac{384865}{19}\omega\xi-\frac{198967}{19}\omega\xi^2.
\]
Since (\ref{rt}) determines $X,Y$ up to a unit factor in $M$, hence
solving the relative Thue equation (\ref{rt}) 
we represent $\beta=X-\xi Y$ using the relative units of $K$ over $M$
and including the torsion units in $\gamma$ (see \cite{gp}, \cite{book}) as 
\[
\beta =\gamma \eta_2^{a_2}\eta_3^{a_3}
\]
where $\gamma=\gamma_i$ with one of $i=1,2,3$ (the calculation must be performed
for all possible $i$) and the exponents $a_2,a_3$ are integers.
Using Baker's method we obtain
\[
A=\max(|a_2|,|a_3|)<10^{29}
\]
which bound is reduced to $A<10$. For all possible values of $a_2,a_3$
we calculated $\gamma \eta_2^{a_2}\eta_3^{a_3}$.
Taking conjugates of the equation
\[
(x_1+\omega y_1)-\xi (x_2+\omega y_2)=\gamma \eta_2^{a_2}\eta_3^{a_3}
\]
we obtain a system of linear equations in $x_1,y_1,x_2,y_2$ which enables
us to calculate the solutions corresponding to $\gamma$, $a_2,a_3$.
(The resolution of this system of linear equations can be
replaced by calculating the representation (\ref{ggg3}) of the product on the
right hand side using integer arithmetic and checking if it
can be written in the form of the left hand side.)

In our example up to sign the only solutions are
\[
(x_{10},y_{10},x_{20},y_{20})=(0,1,0,-1),(0,2,0,1).
\]
Set $X_0=x_{10}+\omega y_{10},Y_0=x_{20}+\omega y_{20}$.
Then all solutions are of the form (\ref{yyyy}).

We consider
\[
I_{K/L}(\alpha)=1,\;\; J_{L,M}(\alpha)=1,
\]
that is
\[
I_{K/L}(\underline{x})=I_{K/L}(x_1,x_2,y_0,y_1,y_2,)=19^{3},\;
J_{L,M}(\underline{x})=J_{L,M}(x_1,x_2,y_0,y_1,y_2)=19^{6}.
\]
The first of these equations is simple:
\[
N_{L/\Q}(y_0+y_1\xi+y_2\xi^2)=\pm 19^2
\]
(we can extract 19 from $I_{K/L}(\underline{x})$).
We set $e=\varepsilon^h$,
substitute the representation (\ref{xy}) of $x_1,y_1,x_2,y_2$ into
the above equations and multiply them by $e^3$ and $e^6$, respectively.
Then we obtain polynomial equations
$F_2(y_0,e)=0$ and $F_3(y_0,e)=0$.
Calculating the resolvent of the above two equations
with respect to $e$ we get a polynomial equation $F_4(y_0)=0$ of degree 72 in $y_0$.
Finding the possible $y_0$, calculating the corresponding $e$ from $F_2(y_0,e)=0$
and calculating $x_1,y_1,x_2,y_2$, we obtain, that up to equivalence the only 
generator of power integral basis in $K$ is 
\[
\alpha=\frac{2\omega+\xi\omega-\xi^2\omega}{19}.
\]

\section{Composites of quartic fields with real quadratic fields}
\label{c42}

Let $M$ be a real quadratic number field with integral basis $\{1,\omega\}$
and fundamental unit $\varepsilon>1$.
Let $L=\Q(\xi)$ be a quartic number field, with a generator element $\xi\in\Z_L$.
In this case the representation (\ref{repr}) is written for simplicity in the
form
\begin{equation}
\alpha=\frac{x_0+x_1\xi+x_2\xi^2+x_3\xi^3+
y_0\omega+y_1\omega \xi+y_2\omega \xi^2+y_3\omega \xi^3}{d},
\label{repr8}
\end{equation}
where $x_i,y_i\in\Z\; (0\leq i\leq 3)$ and $d\in\Z$ is a common denominator.

According to \cite{relquartic} (see also \cite{book})
the calculation of $\alpha\in\Z_K$ with 
\[
I_{K/M}(\alpha)=1
\]
can be reduced to the resolution of certain cubic and quartic 
relative Thue equations.
We write (\ref{repr8}) in the form
\begin{equation}
\alpha=\frac{X_0+X_1\xi+X_2\xi^2+X_3\xi^3}{d},
\label{repr82}
\end{equation}
with $X_i=x_i+\omega y_i\in\Z_M,\; (0\leq i\leq 3)$.

Assume $\xi$ has relative minimal polynomial 
$f(x)=x^4+a_1x^3+a_2x^2+a_3x+a_4\in\Z_M[x]$.
(Note that in our case this is just the absolute minimal polynomial of $\xi$,
i.e. $f(x)\in\Z[x]$.)
Let $i_0=I_{K/M}(\xi)=(\Z_K^+:\Z_M[\xi]^+)$, let
\[
F(u,v)=u^3-a_2u^2v+(a_1a_3-4a_4)uv^2+(4a_2a_4-a_3^2-a_1^2a_4)v^3
\]
be a binary cubic form and let
\begin{eqnarray*}
Q_1(x_1,x_2,x_3)&=&
x_1^2-x_1x_2a_1+x_2^2a_2+x_1x_3(a_1^2-2a_2)
+x_2x_3(a_3-a_1a_2)+x_3^2(-a_1a_3+a_2^2+a_4),\\
Q_2(x_1,x_2,x_3)&=&x_2^2-x_1x_3-a_1x_2x_3+x_3^2a_2,
\end{eqnarray*}
be ternary quadratic forms.

\begin{lemma}(\cite{relquartic})
If $\alpha$ of the form (\ref{repr82}) satisfies
\[
I_{K/M}(\alpha)=1,
\]
then there is a solution $(U,V)\in \Z_M$ of
\begin{equation}
N_{M/\Q}(F(U,V))=\pm \frac{d^{12}}{i_0}
\label{res8}
\end{equation}
such that
\begin{eqnarray}
U&=&Q_1(X_1,X_2,X_3), \nonumber \\
V&=&Q_2(X_1,X_2,X_3).   \label{QV8}
\end{eqnarray}
\label{lllmmm}
\end{lemma}

Equation (\ref{res8}) is either reducible, or it is a cubic
relative Thue equation over $M$. In either case it allows us to determine 
$U,V$ up to a unit factor in $M$:
\begin{equation}
U=\varepsilon^{h_0}\cdot U_{0},\; V=\varepsilon^{h_0}\cdot V_{0}
\label{r123}
\end{equation}
where the finitely many candidates for $U_{0}, V_{0}\in\Z_M$ 
can be calculated.
We set $h_0=2\cdot h+r_0$ with $h\in\Z,r_0\in\{0,1\}$, further let
$X_{i0}=\varepsilon^{-h}X_i\; (1\leq i\leq 3)$.
Then (\ref{QV8}) implies
\begin{eqnarray}
U_0\varepsilon^{r_0}&=&Q_1(X_{10},X_{20},X_{30}), \nonumber \\
V_0\varepsilon^{r_0}&=&Q_2(X_{10},X_{20},X_{30}),   \label{QV8r}
\end{eqnarray}
where the possible values of 
$U_0\varepsilon^{r_0}, V_0\varepsilon^{r_0}$ are known.

Following the arguments of \cite{relquartic} (see also \cite{book}),
some common multiples of $X_{10}, X_{20},X_{30}$, say
$\delta X_{10}, \delta X_{20},\delta X_{30}$ can be represented as 
quadratic forms in 
certain parameters $P,Q\in\Z_M$, where the common factor $\delta$
may only take finitely many values: 
\begin{equation}
\delta X_{10}=f_1(P,Q),\;\;
\delta X_{20}=f_2(P,Q),\;\;
\delta X_{30}=f_3(P,Q).
\label{fpq}
\end{equation}
Substituting these representations of 
$\delta X_{10}, \delta X_{20},\delta X_{30}$ into (\ref{QV8r})
we obtain homogeneous quartic equations in $P,Q$:
\begin{equation}
F_1(P,Q)=\delta^2\ U_0\;\varepsilon^{r_0}, \;\;
F_2(P,Q)=\delta^2\ V_0\;\varepsilon^{r_0}.
\label{r0123}
\end{equation}
According to \cite{relquartic} (see also \cite{book}) at least one of these
equations is a quartic relative Thue equation in $P,Q$
over $M$. The right hand
sides are known, therefore using the method of \cite{gp} all solutions
$P,Q$ can be determined. This enables us to determine $X_{10},X_{20},X_{30}$
using their representation (\ref{fpq}) as quadratic forms in $P,Q$.
This means that we obtain $X_i$ up to a unit factor:
\begin{equation}
X_i=\varepsilon^{h}\ X_{i0}\; (1\leq i\leq 3).
\label{xxx}
\end{equation}

Finally, we have to determine the exponent $h$ of $\varepsilon^h$ and $y_0$ of
$X_0=x_0+\omega y_0$. To do this we proceed similarly as in Section \ref{63}.
(\ref{xxx}) implies
\begin{eqnarray}
x_i&=& \pm 
\frac{\varepsilon^h \overline{\omega}X_{i0}
-\overline{\varepsilon}^h \omega \overline{X_{i0}}}
{\overline{\omega}-\omega},\nonumber \\
y_i&=&\pm
\frac{\varepsilon^h X_{i0}
-\overline{\varepsilon}^h  \overline{X_{i0}}}
{\omega-\overline{\omega}},\label{x3} \\  \nonumber
\end{eqnarray}
for $i=1,2,3$, where overline denotes the conjugate of an element of $M$.

We consider
\[
I_{K/L}(\alpha)=1,\;\; J_{L,M}(\alpha)=1,
\]
that is
\[
I_{K/L}(\underline{X})=I_{K/L}(x_1,x_2,x_3,y_0,y_1,y_2,y_3)=d^{4},
\]
\begin{equation}
J_{L,M}(\underline{X})=J_{L,M}(x_1,x_2,x_3,y_0,y_1,y_2,y_3)=d^{12}.
\label{ijij}
\end{equation}
Set $e=\varepsilon^h$.
If we substitute the expressions (\ref{x3}) of $x_i,y_i\; (1\leq i\leq 3)$
into the above equations
and we multiply them by $e^4$ and $e^{12}$, respectively,
 then we obtain polynomial equations in $e,y_0$:
\begin{equation}
F_2(e,y_0)=0,\; F_3(e,y_0)=0.
\label{ee4}
\end{equation}
Here $F_2$ is of degree 8 in $e$ and of degree 4 in $y_0$, 
$F_3$ is of degree 24 in $e$ and of degree 12 in $y_0$.
Let $F_4(y_0)$ be the resultant of $F_2$ and $F_3$ with respect to $e$.
This is a polynomial equation in $y_0$ with rational coefficients of degree 192. 
It is easy to find the integer roots $y_0$ of $F_4(y_0)=0$.
For a given $y_0$ we calculate the real roots of
\[
F_2(e,y_0)=0
\]
in $e$ and for all real roots $e=\gamma$ we check if
\[
\frac{\ln(|\gamma|)}{\ln\varepsilon}
\]
is an integer. If so, then this value is a candidate for $h$.
We obtain $x_1,y_1,x_2,y_2,x_3,y_3$ from (\ref{x3}).
We have to check if $\alpha$ of (\ref{repr8}) is an algebraic 
integer with some $x_0\in\Z$
and if it has $I(\alpha)=1$.

\section{Composites of totally complex quartic fields \\ with real quadratic fields}

According to \cite{gppjsc} (see also \cite{order}, \cite{book}) the resolution of
index form equations in totally complex quartic fields can be made easier
by using the fact that in that case there is a $\lambda$ such that the 
linear combination $Q_1+\lambda Q_2$ of the quadratic forms $Q_1,Q_2$,
analogous to those in Lemma \ref{lllmmm}, is positive definite. Here we extend this
property to our relative case.

Similarly as above let $f(x)=x^4+a_1x^3+a_2x^2+a_3x+a_4\in\Z[x]$
be the defining polynomial of $\xi$, having four complex roots.
Let $F(u,v),Q_1(x_1,x_2,x_3)$ and $Q_1(x_1,x_2,x_3)$ be the same as
in Lemma \ref{lllmmm}.

By the results of \cite{gppjsc}, $F(u,1)$ 
has three real roots $\lambda_1<\lambda_2<\lambda_3$. If
$-\lambda_2<\lambda<-\lambda_1$, then
$Q_1(x_1,x_2,x_3)+\lambda Q_2(x_1,x_2,x_3)$ is a positive definite quadratic form.

In the relative case we have:

\begin{lemma}
\label{qq}
Let $\lambda_1,\lambda_2,\lambda_3$ be the (real) roots of $F(u,1)$ and
let $-\lambda_2<\lambda<-\lambda_1$.
Set $X_j=x_j+\omega y_j\ (1\leq j\leq 3)$. 
Then $Q_1(X_1,X_2,X_3)+\lambda Q_2(X_1,X_2,X_3)$ can be written in the form
\[
S(x_1,x_2,x_3,y_1,y_2,y_3)+\sqrt{d}\cdot T(x_1,x_2,x_3,y_1,y_2,y_3)
\]
where $S,T$ are quadratic forms with rational coefficients and $S$ is
positive definite.
\end{lemma}

\noindent
{\bf Proof}\\
Let
\[
Q(x_1,x_2,x_3)=\sum_{i=1}^2 \sum_{j=1}^2 a_{ij} x_ix_j
\]
be any quadratic form with integer coefficients. 

First consider the case when $\omega=\sqrt{d}$.

If we substitute $X_i=x_i+\sqrt{d} y_i$, then
\[
Q(X_1,X_2,X_3)=\sum_{i=1}^2 \sum_{j=1}^2 a_{ij} (x_i+\sqrt{d} y_i)(x_j+\sqrt{d} y_j)
\]
\[
=\sum_{i=1}^2 \sum_{j=1}^2 a_{ij} x_ix_j
+ d\sum_{i=1}^2 \sum_{j=1}^2 a_{ij}y_iy_j
+\sqrt{d}\sum_{i=1}^2 \sum_{j=1}^2a_{ij}(x_iy_j+y_ix_j)
\]
\[
=Q(x_1,x_2,x_3)+dQ(y_1,y_2,y_3)+\sqrt{d}\sum_{i=1}^2 \sum_{j=1}^2a_{ij}(x_iy_j+y_ix_j)
\]
Therefore
\[
Q_1(X_1,X_2,X_3)+\lambda Q_2(X_1,X_2,X_3)
\]
\[
=Q_1(x_1,x_2,x_3)+\lambda Q_2(x_1,x_2,x_3)+d[Q_1(y_1,y_2,y_3)+\lambda Q_2(y_1,y_2,y_3)]
\]
\[
+\sqrt{d}\cdot T(x_1,x_2,x_3,y_1,y_2,y_3),
\]
where $T$ is a quadratic form with integer coefficients.
By \cite{gppjsc}, $Q_1(x_1,x_2,x_3)+\lambda Q_2(x_1,x_2,x_3)$
and $Q_1(y_1,y_2,y_3)+\lambda Q_2(y_1,y_2,y_3)$ are positive definite, 
which implies that $S(x_1,x_2,x_3,y_1,y_2,y_3)$ is positive definite, since $d>0$.

Let now $\omega=(1+\sqrt{d})/2$.

If we substitute $X_i=x_i+y_i(1+\sqrt{d})/2$, then
\[
Q(X_1,X_2,X_3)
=\frac{1}{4}Q(2x_1+y_1+\sqrt{d}y_1,2x_2+y_2+\sqrt{d}y_2,2x_3+y_3+\sqrt{d}y_3)
\]
\[
=\frac{1}{4}Q(2x_1+y_1,2x_2+y_2,2x_3+y_3)+\frac{1}{4} d Q(y_1,y_2,y_3)
+\frac{1}{4}\sqrt{d}\sum_{i=1}^2 \sum_{j=1}^2a_{ij}((2x_i+y_i)y_j+y_i(2x_j+y_j)).
\]
We obtain
\[
Q_1(X_1,X_2,X_3)+\lambda Q_2(X_1,X_2,X_3)
\]
\[
=\frac{1}{4}\left[
Q_1(2x_1+y_1,2x_2+y_2,2x_3+y_3)+\lambda Q_2(2x_1+y_1,2x_2+y_2,2x_3+y_3)
\right]
\]
\[
+\frac{1}{4} d \left[
(Q_1(y_1,y_2,y_3)+\lambda Q_2(y_1,y_2,y_3))
\right]
\]
\[
+\sqrt{d}\cdot T(x_1,x_2,x_3,y_1,y_2,y_3).
\]
where $T$ is a quadratic form with rational coefficients
($4T$ has integer coefficients).
By \cite{gppjsc}, $Q_1(2x_1+y_1,2x_2+y_2,2x_3+y_3)
+\lambda Q_2(2x_1+y_1,2x_2+y_2,2x_3+y_3)$
and $Q_1(y_1,y_2,y_3)+\lambda Q_2(y_1,y_2,y_3)$ are positive definite, 
which implies that $S(x_1,x_2,x_3,y_1,y_2,y_3)$ is positive definite, since $d>0$.
$\Box$\\

If the above Lemma \ref{qq} is applicable (i.e. $L$ is totally complex)
then it makes our calculation much easier. 
For given $U_0,V_0$ (cf. (\ref{r123})) 
to find the corresponding $X_{10},X_{20},X_{30}$  (see (\ref{QV8r}))
we do not need to represent $X_{10},X_{20},X_{30}$ 
as quadratic forms in certain parameters $P,Q\in\Z_M$ (cf. (\ref{fpq}))
and then to solve
the corresponding quartic relative Thue equations (\ref{r0123}).
We simply set
\[
Q_1(X_{10},X_{20},X_{30})+\lambda Q_2(X_{10},X_{20},X_{30})
=(U_0+\lambda V_0)\varepsilon^{r_0}
\]
we take the $\sqrt{d}$-free parts of both sides and enumerate 
all possible components $x_{i0}, y_{i0}$ of $X_{i0}=x_{i0}+\omega y_{i0}
\ (1\leq i\leq 3)$.
Having $X_{i0}\ (1\leq i\leq 3)$ we proceed similarly as in Section \ref{c42}
determining $y_0$ and $\varepsilon^h$.

\section{Example 2}

Let $M=\Q(\sqrt{2})$ with fundamental unit $\varepsilon=(1+\sqrt{2})$.
Set $\omega=\sqrt{2}$. Let $L=\Q(\xi)$ where $\xi$ has minimal polynomial
$x^4+2x^2+2x+1$. 
The composite field $K=LM=\Q(\xi,\sqrt{2})$ is generated
by $\vartheta=\xi\sqrt{2}$ 
with minimal polynomial $x^8+8x^6+24x^4+16$.
Considering the integral basis of $K$ we may represent all $\alpha\in\Z_K$
in the from
\[
\alpha=\frac{x_0+x_1\xi+x_2\xi^2+x_3\xi^3+
y_0\omega+y_1\omega \xi+y_2\omega \xi^2+y_3\omega \xi^3}{4}.
\]
Our purpose is to determine
all $\alpha\in\Z_K$ having index 1, that is, generating a power integral basis.

We have 
\[
D_M=2^3,\; D_L=2^4\cdot 37,\; D_K=2^{16}\cdot 37^2.
\]
By 
\[
D_K=N_{M/\Q}(D_{K/M})\cdot D_M^4,\;\; D_K=N_{L/\Q}(D_{K/L})\cdot D_L^2
\]
we obtain
\[
N_{M/\Q}(D_{K/M})=2^4\cdot 37^2, \;\; N_{L/\Q}(D_{K/L})=2^8.
\]
The relative Thue equation (\ref{res8}) has the form
\[
N_{M/\Q}(F(U,V))=\pm 2^{22}.
\]
with 
\[
F(U,V)=U^3-2U^2V-4UV^2+4V^3.
\]
Let $\rho$ be a root of $F(x,1)$. 
The element  $\vartheta=\rho\sqrt{2}$
has minimal polynomial $x^6-24x^4+128x^2-128$. 
Set  $H=\Q(\rho\sqrt{2})$.
In $\Z_H$ up to associates there is one element of norm $2^{22}$:
\[
\gamma=-8\rho.
\]

The fundamental units of the sextic field $H$  are
\begin{eqnarray*}
\eta_1&=&\varepsilon=1+\frac{4\vartheta}{2}-\frac{5\vartheta^3}{8}+
\frac{\vartheta^5}{32},\\
\eta_2&=&-1+\frac{\vartheta^2}{4},\\
\eta_3&=&-1+\frac{3\vartheta}{2}-\frac{5\vartheta^3}{8}+
\frac{\vartheta^5}{32},\\
\eta_4&=&1-\frac{\vartheta}{2},\\
\eta_5&=&1+\frac{3\vartheta}{2}-\frac{5\vartheta^3}{8}+
\frac{\vartheta^5}{32}.
\end{eqnarray*}
Solving the relative Thue equation (\ref{res8}) we represent $\beta=U-\rho V$
as 
\[
\beta =\gamma \eta_2^{a_2}\eta_3^{a_3}\eta_4^{a_4}\eta_5^{a_5}
\]
where the exponents $a_2,a_3,a_4,a_5$ are integers
(cf. the corresponding remarks in Example 1).
Using Baker's method we obtain
\[
A=\max(|a_2|,|a_3|,|a_4|,|a_5|)<10^{33}
\]
which bound is reduced to $A<22$. For all possible values of $a_2,a_3,a_4,a_5$
we calculated $\gamma \eta_2^{a_2}\eta_3^{a_3}\eta_4^{a_4}\eta_5^{a_5}$.
Taking conjugates of the equation
\[
(u_{10}+\omega u_{20})-\rho (v_{10}+\omega v_{20})
=\gamma \eta_2^{a_2}\eta_3^{a_3}\eta_4^{a_4}\eta_5^{a_5}
\]
we obtain a system of linear equations in $u_{10},u_{20},v_{10},v_{20}$ which enables
us to calculate the solutions corresponding to $\gamma$, $a_2,a_3,a_4,a_5$.
(Similarly to Example 1, this can also be done by using only integer arithmetic.)
Up to sign we obtained the following solutions:
\[
(u_{10},u_{20},v_{10},v_{20})=
(720,0,248,0),(32,16,8,8),(32,-16,8,-8),
(16,0,8,0),(16,0,-8,0),
\]
\[(0,0,8,0),(16,0,8,8),(64,-16,40,8),(16,0,8,-8),
(64,16,40,-8),(32,0,40,0).
\]
The roots of $F(u,1)$ are approximately
$\lambda_1=-1.709$,
$\lambda_2=0.806$ and
$\lambda_3=2.903$.
A suitable value of $\lambda$ is $\lambda=0$.
For all the above possible $U_0=u_{10}+\omega u_{20},V_0=v_{10}+\omega v_{20}$
we consider 
\[
Q_1(X_{10},X_{20},X_{30})+\lambda Q_2(X_{10},X_{20},X_{30})=
(U_0+\lambda V_0)\varepsilon^{r_0},
\]
i.e.
\[
Q_1(X_{10},X_{20},X_{30})=U_0\varepsilon^{r_0},
\]
with $r_0\in\{0,1\}$.
We take the $\sqrt{2}$-free parts on both sides and using Lemma \ref{qq} 
we enumerate all possible $x_{10},x_{20},x_{30},y_{10},y_{20},y_{30}$.
We set $X_{i0}=x_{i0}+\omega y_{i0}$ $(1\leq i\leq 3)$. 
By $X_i=x_i+\omega y_i
=\varepsilon^h X_{i0}$
we obtain $x_i,y_i$ in the representation (\ref{x3}).

We consider
\[
I_{K/L}(\alpha)=1,\;\; J_{L,M}(\alpha)=1,
\]
that is
\[
I_{K/L}(\underline{x})=I_{K/L}(x_1,x_2,x_3,y_0,y_1,y_2,y_3)=4^{4},\;
J_{L,M}(\underline{x})=J_{L,M}(x_1,x_2,x_3,y_0,y_1,y_2,y_3)=4^{12}.
\]

The first of these equations is simple:
\[
N_{L/\Q}(y_0+y_1\xi+y_2\xi^2++y_3\xi^3)=\pm 2^6
\]
of degree 4 in $y_0$ (a factor 4 can be extracted from $I_{K/L}(\underline{x})$)
and the second equation is of degree 12 in $y_0$.
We set $e=\varepsilon^h$,
substitute the representation (\ref{x3}) of $x_1,y_1,x_2,y_2,x_3,y_3$ into
the above equations and multiply them by $e^4$ and $e^{12}$, respectively.
Then we obtain polynomial equations
$F_2(y_0,e)=0$ and $F_3(y_0,e)=0$.
Calculating the resolvent of the above two equations
with respect to $e$ we get a polynomial equation of degree 192 in $y_0$.
Finding the possible $y_0$, and calculating the corresponding $e$ from $F_2(y_0,e)=0$
and $x_1,y_1,x_2,y_2,x_3,y_3$ we obtain, 
that up to equivalence there are two
generators of power integral basis in $K$:
\[
\alpha=\frac{2\omega+2\omega\xi^2}{4},\; \frac{4\omega+2\omega\xi+2\omega\xi^3}{4}.
\]

\section{Computational aspects}

The field data to our examples were calculated by using KANT \cite{kant}
and Magma \cite{magma}, all other calculations were performed by Maple 
\cite{maple}
on an average Windows 8 laptop with Intel Core i5-4200U CPU having
clock speed 1.3GHz and 4Gb RAM. The calculations took just a couple of minutes.
We especially took care (using high precision)
at determining the real roots $e$ of a relatively high degree
polynomial in the last part of our method. \\

\noindent
{\bf Remark 1} These ideas can easily be 
extended to determining elements of given index
(or minimal index) in these types of number fields.\\

\noindent
{\bf Remark 2}
The method can be easily adopted to the case when 
instead of $\Z_K$ we consider an order $\O$ of $\Z_K$.


\begin{thebibliography}{10}

\normalsize
\baselineskip=17pt

\bibitem{baker}
A. Baker, Transcendental Number Theory,
Cambridge, 1990.


\bibitem{s6} Y. Bilu, I. Ga\'al and K. Gy\H ory, 
{\em Index form equations in sextic 
fields: a hard computation}, Acta Arith., {\bf 115} (2004), No. 1, 85-96.

\bibitem{magma}
W. Bosma, J. Cannon and C. Playoust, {\em The Magma algebra system.
I. The user language}, J. Symbolic Comput., {\bf 24} (1997), 235-265.


\bibitem{maple}
B. W. Char, K. O. Geddes, G. H. Gonnet, M. B. Monagan, S. M. Watt (eds.),
MAPLE, Reference Manual, Watcom Publications, Waterloo, Canada, 1988.

\bibitem{kant}
M. Daberkow, C. Fieker, J. Kl\"uners, M. Pohst, K. Roegner and K. Wildanger,
KANT V4, J. Symbolic Comp., {\bf 24} (1997), 267--283.

\bibitem{fj18}
Z.Franu\u si\' c and B.Jadrijevi\'c,
{\em Computing relative power integral bases in a family 
of quartic extensions of imaginary quadratic fields},
Publ. Math. (Debrecen), {\bf 92}(2018), No. 3-4, 293-315.

\bibitem{order}
I. Ga\'al, {\em Computing all power integral bases in orders of
totally real cyclic sextic number fields},
Math. Comp., {\bf 65} (1996), 801--822.

\bibitem{compos} 
I. Ga\'al, {\em Power integral bases in composites of number fields},
Canad. Math. Bull., {\bf 41} (1998), 158--161.

\bibitem{degnine} 
I. Ga\'al, {\em Solving index form equations in fields of degree nine
with cubic subfields}, J. Symbolic Comput., {\bf 30} (2000), 181-193.

\bibitem{relcubic}
I. Ga\'al, {\em Power integral bases in cubic relative extensions},
Experimental Math., {\bf 10} (2001), 133--139.


\bibitem{book}
I. Ga\'al,
Diophantine equations and power integral bases,
Boston, Birkh\"auser, 2002.

\bibitem{relthuesmall}
I. Ga\'al, {\em Calculating "small" solutions of relative Thue equations}
Experimental Math., {\bf 24} (2015), 142-149.

\bibitem{s5}
I. Ga\'al and K. Gy\H{o}ry, {\em Index form equations in quintic fields},
Acta Arith., {\bf 89} (1999), 379--396.

\bibitem{gop}
I. Ga\'al, P. Olajos and M. Pohst, {\em Power integral bases
in orders of composites of number fields}, Experimental Math.,
{\bf 11} (2002), 87--90.


\bibitem{gppjsc}
I. Ga\'{a}l, A. Peth\H{o} and M. Pohst, {\em On the resolution of index form
equations in quartic number fields}, J. Symbolic Comput.,
{\bf 16} (1993), 563--584.


\bibitem{gppsim}
I. Ga\'{a}l, A. Peth\H{o} and M. Pohst, {\em Simultaneous representation
of integers by a pair of ternary quadratic forms -- with an application to
index form equations in quartic number fields}, J. Number Theory,
{\bf 57} (1996), 90--104.


\bibitem{gpezaz} 
I. Ga\'{a}l and M. Pohst,  {\em On the resolution of index form
equations in sextic fields with an imaginary quadratic subfield},
J. Symbolic Comp., {\bf 22} (1996), 425--434.

\bibitem{relquartic} 
I. Ga\'al and M. Pohst, {\em On the resolution of index form equations 
in relative quartic extensions}, J. Number Theory, {\bf 85} (2000), 201--219.


\bibitem{gp}
I. Ga\'al and M. Pohst,
{\em On the resolution of relative Thue equations},
Math. Comput., {\bf 71} (2002), 429-440.

\bibitem{s62}
I. Ga\'al and L. Remete,
{\em Power integral bases in a family of sextic fields with quadratic subfields},
Tatra Mt. Math. Publ., {\bf 64} (2015), 59-66.

\bibitem{composite}
I. Ga\'al and L. Remete, {\em Integral bases and monogenity of composite fields},
Experimental Math., online appeared 2017.

\bibitem{simlestsextic}
I. Ga\'al and L. Remete, {\em Integral bases and monogenity of the
simplest sextic fields}, Acta Arith., {\bf 183} (2018), No. 2, 173-183

\bibitem{jalfa}
I. Ga\'al, L. Remete and T. Szab\'o,
{\em Calculating power integral bases by using relative power integral bases}
Functiones et Approximatio Comment. Math., {\bf 54} (2016), No. 2., 141-149.

\bibitem{gsch}
I. Ga\'al and N. Schulte, {\em Computing all power integral bases of
cubic number fields}, Math. Comput., {\bf 53} (1989), 689--696.

\bibitem{gsz2}
I. Ga\'al  and T. Szab\'o, {\em Relative power integral bases
in infinite families of quartic extensions of quadratic field}, 
JP Journal of Algebra, Number Theory and Appl., {\bf 29} (2013), 31--43.

\bibitem{olaj}
P. Olajos, {\it Power integral bases in orders of composite fields. II.},
Ann. Univ. Sci. Budap. Rolando Eötvös, Sect. Math., {\bf 46} (2003), 35-41.

\bibitem{thue}
A. Thue, {\em \"Uber Ann\"aherungswerte algebraischer Zahlen},
J.Reine Angew. Math., {\bf 135} (1909), 284--305.



\end{thebibliography}
\end{document}